\input amstex
\input amsppt.sty
\magnification=\magstep1
\hsize=32truecc
\vsize=22.5truecm
\baselineskip=16truept
\NoBlackBoxes
\TagsOnRight \pageno=1 \nologo
\def\Z{\Bbb Z}

\def\l{\left}
\def\r{\right}
\def\bg{\bigg}
\def\({\bg(}
\def\[{\bg\lfloor}
\def\){\bg)}
\def\]{\bg\rfloor}
\def\t{\text}
\def\f{\frac}

\def\ls{\leqslant}
\def\gs{\geqslant}

\def\al{\alpha}

\def\Proof{\noindent{\it Proof}}

\def\Remark{\medskip\noindent{\it  Remark}}

\def\Ack{\medskip\noindent {\bf Acknowledgment}}
\hbox {Bull. Aust. Math. Soc. 88(2013), no.\,2, 197--205.}
\bigskip
\topmatter
\title On a sequence involving sums of primes\endtitle
\author Zhi-Wei Sun\endauthor
\leftheadtext{Zhi-Wei Sun}
 \rightheadtext{On a sequence involving sums of primes}
\affil Department of Mathematics, Nanjing University\\
 Nanjing 210093, People's Republic of China
  \\  zwsun\@nju.edu.cn
  \\ {\tt http://math.nju.edu.cn/$\sim$zwsun}
\endaffil
\abstract For $n=1,2,3,\ldots$ let $S_n$ be the sum of the first $n$
primes. We mainly show that the sequence $a_n=\root n\of{S_n/n}\ (n=1,2,3,\ldots)$ is strictly decreasing, and moreover the
sequence $a_{n+1}/a_n\ (n=10,11,\ldots)$ is strictly increasing.
We also formulate similar conjectures involving twin primes or partitions of integers.
\endabstract
\thanks 2010 {\it Mathematics Subject Classification}.\,Primary 11A41;
Secondary 05A17, 05A20, 11B75, 11B83, 11J99, 11P83.
\newline\indent Supported by the National Natural Science
Foundation (grant 11171140) of China and the PAPD of Jiangsu Higher
Education Institutions.
\endthanks
\endtopmatter
\document

\heading{1. Introduction}\endheading

For $n\in\Z^+=\{1,2,3,\ldots\}$ let $p_n$ denote the $n$th prime.
The unsolved Firoozbakht conjecture (cf. [R, p.\,185]) asserts that
$$\root n\of{p_n}>\root{n+1}\of {p_{n+1}}\quad\t{for all}\ n\in\Z^+,$$
i.e., the sequence $(\root{n}\of {p_n})_{n\gs1}$ is strictly decreasing.
This implies the inequality $p_{n+1}-p_n<\log^2 p_n-\log p_n+1$ for large $n$, which is even stronger than Cram\'er's conjecture
$p_{n+1}-p_n=O(\log^2 p_n)$.
Let $P_n$ be the product of the first $n$ primes. Then $P_n<p_{n+1}^n$ and hence
$P_n^{n+1}<P_{n+1}^n$. So the sequence $(\root{n}\of{P_n})_{n\gs1}$ is strictly increasing.

Now let us look at a simple example not related to primes.

{\it Example}\ 1.1. Let $a_n=\root n\of{n}$ for $n\in\Z^+$. Then the sequence
$(a_n)_{n\gs 3}$ is strictly decreasing, and the sequence $(a_{n+1}/a_n)_{n\gs4}$
is strictly increasing. To see this we investigate the function $f(x)=\log(x^{1/x})=(\log x)/x$ with $x\gs3$.
As $f'(x)=(1-\log x)/x^2<0$, we have $f(n)>f(n+1)$ for $n=3,4,\ldots$. Since
$$f''(x)=\f{2\log x-3}{x^3}>0\quad\t{for}\ x\gs 4.5,$$
the function $f(x)$ is strictly convex over the interval $(4.5,+\infty)$ and so
$$2f(n+1)<f(n)+f(n+2)\ (\t{i.e.,}\ a_{n+1}^2<a_na_{n+2})\quad \t{for}\ n=5,6,\ldots.$$
The inequality $a_5^2<a_4a_6$ can be verified directly.
\medskip

A sequence $(a_n)_{n\gs1}$ of nonnegative real numbers is said to be {\it log-convex} if $a_{n+1}^2\ls a_na_{n+2}$ for all $n=1,2,3,\ldots$.
Many combinatorial sequences (such as the sequence of the Catalan numbers) are log-convex,
the reader may consult [LW] for some results on log-convex sequences.

For $n\in \Z^+$ let $S_n=\sum_{k=1}^np_k$ be the
sum of the first $n$ primes.
For instance,
$$S_1=2,\ S_2=2+3=5,\ S_3=2+3+5=10,\ S_4=2+3+5+7=17.$$
Recently the author [S] conjectured that for any positive integer
$n$ the interval $(S_n,S_{n+1})$ contains a prime. As $S_n<np_{n+1}$ for
all $n\in\Z^+$,  the sequence $(S_n/n)_{n\gs1}$ is strictly
increasing.

In the next section we will state our theorems involving the sequence $(a_n)_{n\gs1}$ with
$a_n=\root n\of{S_n/n}$, and pose three related conjectures for further research.
Section 3 is devoted to our proofs of the theorems.

\heading{2. Our results and conjectures}\endheading

\proclaim{Theorem 2.1} The sequences $(\root n\of{S_n})_{n\gs2}$ and $(\root n\of{S_n/n})_{n\gs1}$ are strictly decreasing.
\endproclaim
\Remark\ 2.2. Note that $S_n/n$ is just the arithmetic mean of the first $n$ primes.
It is interesting to compare Theorem 2.1 with Firoozbakht's conjecture that $(\root n\of{p_n})_{n\gs1}$ is strictly decreasing.
\medskip

For $\al>0$ and $n\in\Z^+$ define
$$S_n^{(\al)}=\sum_{k=1}^np_k^{\al}.$$

We actually obtain the following extension of Theorem 2.1.

\proclaim{Theorem 2.3} Let $\al\gs1$ and
$n\in\Z^+$ with $n\gs\max\{100,e^{2\times1.348^{\al}+1}\}$. Then
$$\root n\of{\f{S_n^{(\al)}}n}>\root{n+1}\of{\f{S_{n+1}^{(\al)}}{n+1}}\tag2.1$$
and hence
$$\root n\of{S_n^{(\al)}}>\root{n+1}\of{S_{n+1}^{(\al)}}.\tag2.2$$
\endproclaim
\Remark\ 2.4. In view of Example 1.1, (2.1) implies (2.2) if $n\gs3$. We conjecture that (2.1) holds for any $\al>0$ and $n\in\Z^+$.
\medskip

Note that $\lfloor e^{2\times 1.348+1}\rfloor=40$ and we can easily verify that
$$\root{n}\of{\f{S_n}n}>\root{n+1}\of{\f{S_{n+1}}{n+1}}\quad\t{for every}\ n=1,\ldots,99.$$
So Theorem 2.1 follows from Theorem 2.3 in the case $\al=1$.

\proclaim{Corollary 2.5} For each $\al\in\{2,3,4\}$, the sequences
$$\l(\root n\of{\f{S_n^{(\al)}}n}\,\r)_{n\gs1}\ \ \t{and}\ \ \ \l(\root n\of{S_n^{(\al)}}\,\r)_{n\gs1}$$
are strictly decreasing.
\endproclaim
\Proof. Observe that
$$\lfloor e^{2\times 1.348^2+1}\rfloor=102,\ \lfloor e^{2\times 1.348^3+1}\rfloor=364,\ \lfloor e^{2\times 1.348^4+1}\rfloor=2005.$$
In light of Theorem 2.3 and Example 1.1, it suffices to verify that
$$\root n\of{\f{S_n^{(\al)}}n}>\root {n+1}\of{\f{S_{n+1}^{(2)}}{n+1}}$$
whenever $\al\in\{2,3,4\}$ and $n\in\{1,\ldots,\lfloor e^{2\times 1.348^\al+1}\rfloor\}.$
These can be easily done via computer. \qed
\medskip

Our following theorem is more sophisticated than Theorem 2.3.

\proclaim{Theorem 2.6} Let $\al\gs1$.
Then the sequence
$$\l(\root{n+1}\of{S^{(\al)}_{n+1}/(n+1)}\bigg/\root{n}\of{S^{(\al)}_n/n}\r)_{n\gs N(\al)}$$ is strictly increasing,
where
$$N(\al)=\max\l\{350000,\ \lceil e^{((\al+1)^21.2^{2\al+1}+(\al+1)1.2^{\al+1})/\al}\rceil\r\}.\tag2.3$$
\endproclaim

\proclaim{Corollary 2.7} All the sequences
$$\gather\l(\root{n+1}\of{S_{n+1}/(n+1)}\big/\root{n}\of{S_n/n}\r)_{n\gs10},\ \l(\root{n+1}\of{S_{n+1}}\big/\root{n}\of{S_n}\r)_{n\gs5},
\\ \l(\root{n+1}\of{S^{(2)}_{n+1}/(n+1)}\big/\root{n}\of{S^{(2)}_n/n}\r)_{n\gs13},
\ \l(\root{n+1}\of{S^{(2)}_{n+1}}\big/\root{n}\of{S^{(2)}_n}\r)_{n\gs10},
\\ \l(\root{n+1}\of{S^{(3)}_{n+1}/(n+1)}\big/\root{n}\of{S^{(3)}_n/n}\r)_{n\gs17},
\ \l(\root{n+1}\of{S^{(3)}_{n+1}}\big/\root{n}\of{S^{(3)}_n}\r)_{n\gs10},
\\ \l(\root{n+1}\of{S^{(4)}_{n+1}/(n+1)}\big/\root{n}\of{S^{(4)}_n/n}\r)_{n\gs35},
\ \l(\root{n+1}\of{S^{(4)}_{n+1}}\big/\root{n}\of{S^{(4)}_n}\r)_{n\gs17}
\endgather$$
are strictly increasing.
\endproclaim
\Proof. For $N(\al)$ given by (2.3), via computation we find that
$$N(1)=350000,\ N(2)=974267,\ N(3)=3163983273$$
and
$$N(4)=2271069361863763.$$
Via computer we can verify that
$$\f{\root{n+1}\of{S_{n+1}^{(\al)}/(n+1)}}{\root n\of{S_n^{(\al)}/n}}<\f{\root{n+2}\of{S_{n+2}^{(\al)}/(n+2)}}{\root {n+1}\of{S_{n+1}^{(\al)}/(n+1)}}$$
for all $\al\in\{1,2,3,4\}$ and $n=N_0(\al),\ldots,N(\al)-1$, where
$$N_0(1)=10,\ N_0(2)=13,\ N_0(3)=17,\ N_0(4)=35.$$
Combining this with Theorem 2.3 we obtain that $$\l(\root{n+1}\of{S_{n+1}/(n+1)}/\root{n}\of{S_n/n}\r)_{n\gs N_0(\al)}$$ is strictly increasing for each
$\al=1,2,3,4$. Recall that  $(\root{n+1}\of{n+1}/\root n\of{n})_{n\gs 4}$ is strictly increasing by Example 1.1. So
$\l(\root{n+1}\of{S_{n+1}}/\root{n}\of{S_n}\r)_{n\gs N_0(\al)}$ is strictly increasing for any
$\al\in\{1,2,3,4\}$. It remains to check that
$$\f{\root{n+1}\of{S_{n+1}^{(\al)}}}{\root n\of{S_n^{(\al)}}}<\f{\root{n+2}\of{S_{n+2}^{(\al)}}}{\root {n+1}\of{S_{n+1}^{(\al)}}}$$
for all $\al\in\{1,2,3,4\}$ and $n=n_0(\al),\ldots,N_0(\al)-1$, where
$n_0(1)=5$, $n_0(2)=n_0(3)=10$, and $n_0(4)=17$. This can be easily done via computer. \qed
\medskip

We conclude this section by posing three conjectures.

\proclaim{Conjecture 2.8} The two constants
$$s_1=\sum_{n=1}^\infty\f1{S_n}\ \ \t{and}\ \ s_2=\sum_{n=1}^\infty\f{(-1)^n}{S_n}$$
are both transcendental numbers.
\endproclaim
\Remark\ 2.9. Our computation shows that $s_1\approx 1.023476$ and $s_2\approx -0.3624545778$.
\medskip

If $p$ and $p+2$ are both primes, then they are called twin primes.
The famous twin prime conjecture states that there are infinitely many twin primes.

\proclaim{Conjecture 2.10} {\rm (i)} If $\{t_1,t_1+2\},\ldots,\{t_n,t_n+2\}$ are the first $n$ pairs of twin primes, then
the first prime $t_{n+1}$ in the next pair of twin primes is smaller than $t_n^{1+1/n}$, i.e.,
$\root n\of{t_n}>\root{n+1}\of{t_{n+1}}$.

{\rm (ii)} The sequence $(\root{n+1}\of {T_{n+1}}/\root n\of{T_n})_{n\gs9}$ is strictly increasing with limit $1$, where $T_n=\sum_{k=1}^nt_k$.
\endproclaim
\Remark\ 2.11.  Via {\tt Mathematica} the author has verified that $\root n\of{t_n}>\root{n+1}\of{t_{n+1}}$ for all $n=1,\ldots,500000$, and
$\root{n+1}\of {T_{n+1}}/\root n\of{T_n}<\root{n+2}\of {T_{n+2}}/\root {n+1}\of{T_{n+1}}$ for all $n=9,\ldots,500000$. Note that
$t_{500000}=115438667$.
\medskip

Recall that a  partition of a positive integer $n$ is a way of writing $n$ as a sum of  positive integers with the order of addends ignored.
Also, a {\it strict partition} of $n\in\Z^+$ is a way of writing $n$ as a sum of {\it distinct} positive integers with the order of addends ignored.
For $n=1,2,3,\ldots$ we denote by $p(n)$ and $p_*(n)$ the number of partitions of $n$ and the number of strict partitions of $n$ respectively.
It is known that
$$p(n)\sim\f{e^{\pi\sqrt{2n/3}}}{4\sqrt3 n}\ \ \t{and}\ \ p_*(n)\sim\f{e^{\pi\sqrt{n/3}}}{4(3n^3)^{1/4}}\ \quad\t{as}\ n\to+\infty$$
(cf. [HR] and [AS, p.\,826]) and hence $\lim_{n\to\infty}\root{n}\of{p(n)}=\lim_{n\to\infty}\root{n}\of{p_*(n)}=1$.
Here we formulate a conjecture similar to Conjecture 2.10.

\proclaim{Conjecture 2.12} For $n\in\Z^+$ let
$$q(n)=\f{p(n)}n,\ q_*(n)=\f{p_*(n)}n,\ r(n)=\root n\of{q(n)},\ \t{and}\ r_*(n)=\root n\of{q_*(n)}.$$ Then the sequences
$(q(n+1)/q(n))_{n\gs31}$ and $(q_*(n+1)/q_*(n))_{n\gs 44}$ are strictly decreasing, and the sequences
 $(r(n+1)/r(n))_{n\gs60}$ and $(r_*(n+1)/r_*(n))_{n\gs120}$ are strictly increasing.
\endproclaim
\Remark\ 2.13. Via {\tt Mathematica} we have verified the conjecture for $n$ up to $10^5$.
In light of Example 1.1, Conjecture 2.12 implies that all the sequences
$$\l(\f{p(n+1)}{p(n)}\r)_{n\gs25}, \ \l(\f{p_*(n+1)}{p_*(n)}\r)_{n\gs32},\ (\root n\of{p(n)})_{n\gs6},\ (\root n\of{p_*(n)})_{n\gs9}$$
are strictly decreasing, and that the sequences $(\root{n+1}\of{p(n+1)}/\root n\of{p(n)})_{n\gs26}$ and
$(\root{n+1}\of{p_*(n+1)}/\root n\of{p_*(n)})_{n\gs45}$ are strictly increasing. The fact that $(p(n+1)/p(n))_{n\gs 25}$ is strictly decreasing
was conjectured by W.Y.C. Chen [C] and proved by J.E. Janoski [J, pp.\,7-23].

\heading{3. Proofs of Theorems 2.3 and 2.6}\endheading

\proclaim{Lemma 3.1} Let $\al\gs1$ and $n\in\{2,3,\ldots\}$. Then
$$S_n^{(\al)}> 2^\al+\f{n^{\al+1}\log^{\al}n}{\al+1}\l(1-\f{\al}{(\al+1)\log n}\r).\tag3.1$$
\endproclaim
\Proof. It is known that $p_k\gs k\log k$ for $k=2,3,\ldots$ (cf. [Ro] and [RS, (3.12)]). Thus
$$\align S_n^{(\al)}-2^{\al}=\sum_{k=2}^np_k^{\al}\gs&\sum_{k=2}^n(k\log k)^{\al}
>\sum_{k=2}^n\int_{k-1}^k(x\log x)^{\al}dx=\int_1^n(x\log x)^{\al}dx
\endalign$$
Using integration by parts, we find that
$$\align\int_1^n(x\log x)^{\al}dx=&\f{x^{\al+1}}{\al+1}\log^{\al}x\bigg|_{x=1}^n-\int_1^n\l(\f{x^{\al+1}}{\al+1}\cdot\f{\al(\log x)^{\al-1}}x\r)dx
\\=&\f{n^{\al+1}}{\al+1}\log^{\al}n-\f{\al}{\al+1}\int_1^nx^{\al}(\log x)^{\al-1}dx
\\\gs&\f{n^{\al+1}}{\al+1}\log^{\al}n-\f{\al}{\al+1}\int_1^n x^{\al}(\log n)^{\al-1}dx
\\\gs&\f{n^{\al+1}}{\al+1}\log^{\al}n-\f{\al n^{\al+1}}{(\al+1)^2}(\log n)^{\al-1}.
\endalign$$
Therefore (3.1) holds. \qed

\proclaim{Lemma 3.2} Let $\al\gs1$ and $n\in\Z^+$ with $n\gs55$. Then
$$\log S^{(\al)}_n>(\al+1)\log n.\tag3.2$$
\endproclaim
\Proof. Note that $54<e^4<55\ls n$. As $\log^{\al}n>4^{\al}=(2^{\al})^2\gs(\al+1)^2$, by Lemma 3.1 we have
$$S_n^{(\al)}>\f{n^{\al+1}\log^{\al}n}{\al+1}\l(1-\f{\al}{\al+1}\r)=\f{n^{\al+1}}{(\al+1)^2}\log^{\al}n\gs n^{\al+1}$$
and hence (3.2) follows. \qed

\medskip\noindent{\it Proof of Theorem 2.3}.
It is known that
$$p_m<m(\log m+\log\log m)$$
for any $m\gs6$ (cf. [RS, (3.13)] and [D, Lemma 1]). If $m\gs 101$, then
$$\f{\log\log m}{\log m}\ls\f{\log\log 101}{\log101}<0.3314$$
and hence $p_m<1.3314m\log m$. As $n+1\ls 1.01n$, we have
$$\f{\log(n+1)}{\log n}=1+\f{\log((n+1)/n)}{\log n}\ls 1+\f{\log1.01}{\log n}\ls1+\f{\log1.01}{\log100}<1.0022.$$
Therefore
$$p_{n+1}<1.3314(n+1)\log(n+1)<1.3314\times1.01n\times 1.0022\log n<1.348n\log n.$$

Combining Lemmas 3.1 and 3.2, we see that
$$\align &S_n^{(\al)}\l(\f{n+1}{n^{1+1/n}}\root n\of{S_n^{(\al)}}-1\r)
\\=&S_n^{(\al)}\l(e^{(\log S_n^{(\al)})/n+\log(n+1)-(1+1/n)\log n}-1\r)
\\\gs& S_n^{(\al)}\l(e^{(\log S_n^{(\al)}-\log n)/n}-1\r)\gs S_n^{(\al)}\l(e^{(\al\log n)/n}-1\r)
\\>&\f{n^{\al+1}\log^{\al}n}{\al+1}\l(1-\f{\al}{(\al+1)\log n}\r)\f{\al\log n}n
\\=&\f{\al}{\al+1}(n\log n)^{\al}\l(\log n-\f{\al}{\al+1}\r)
\\>&\f{(n\log n)^{\al}}2(\log n-1).
\endalign$$
As $(\log n-1)/2\gs1.348^{\al}$, from the above we get
$$(n+1)\l(\f{S_n^{(\al)}}n\r)^{1+1/n}-S_n^{(\al)}>(1.348n\log n)^{\al}>p_{n+1}^{\al}$$
and hence
$$\l(\f{S_n^{(\al)}}n\r)^{(n+1)/n}>\f{S_{n+1}^{(\al)}}{n+1}$$
which yields (2.1). As mentioned in Remark 2.4, (2.2) follows from (2.1). This concludes the proof. \qed

\medskip\noindent{\it Proof of Theorem 2.6}. Fix an integer $n\gs N(\al)$.
For any integer $m\gs 350001$, we have
$$\f{\log\log m}{\log m}\ls\f{\log\log 350001}{\log 350001}<0.1996$$
and hence
$$p_m<m(\log m)\l(1+\f{\log\log m}{\log m}\r)<1.1996m\log m.$$
 As $n\gs350000$, we have
$$\f{\log(n+1)}{\log n}=1+\f{\log (1+1/n)}{\log n}\ls \f{\log350001}{\log 350000}<1+10^{-6}.$$
Therefore
$$\align p_{n+1}<&1.1996(n+1)\log(n+1)
\\<&1.1996\times\f{350001}{350000}n\times (1+10^{-6})\log n<1.2n\log n.\endalign$$
Since $\log n\gs\log 350000>1/0.078335$, Lemma 3.1 implies that
$$S_n^{(\al)}>\f{n^{\al+1}\log^{\al}n}{\al+1}(1-0.078335)>\f{n^{\al+1}\log^{\al}n}{1.085(\al+1)}.$$
Therefore
$$q_n^{(\al)}:=\f{p_{n+1}^{\al}}{S_n^{(\al)}}<\f{c_{\al}}n,\tag3.3$$
where $c_{\al}=1.085(\al+1)1.2^{\al}$.

By calculus,
$$x-\f{x^2}2<\log(1+x)<x\quad\t{for}\ x>0$$
and
$$-x-x^2<\log(1-x)<-x\ \quad\t{for}\ 0<x<0.5.$$
Thus
$$\log\f{S_{n+1}^{(\al)}/(n+1)}{S_n^{(\al)}/n}=\log\l(1-\f1{n+1}\r)+\log(1+q_n^{(\al)})<-\f1{n+1}+q_n^{(\al)}$$
and
$$\align \log\f{S_{n+2}^{(\al)}/(n+2)}{S_n^{(\al)}/n}>&\log\l(1-\f2{n+2}\r)+\log(1+2q_n^{(\al)})
\\>&-\f2{n+2}-\f4{(n+2)^2}+2q_n^{(\al)}-2(q_n^{(\al)})^2.\endalign$$
Hence
$$\align&D_n^{(\al)}:=\f2{n+1}\log\f{S_{n+1}^{(\al)}}{n+1}-\f1n\log\f{S_n^{(\al)}}n-\f1{n+2}\log\f{S_{n+2}^{(\al)}}{n+2}
\\<&\f 2{n+1}\l(\log\f{S^{(\al)}_n}n-\f1{n+1}+q_n^{(\al)}\r)-\f1n\log\f{S_n^{(\al)}}n
\\&-\f1{n+2}\l(\log\f{S_n^{(\al)}}n-\f2{n+2}-\f4{(n+2)^2}+2q_n^{(\al)}-2(q_n^{(\al)})^2\r)
\\=&\f{-2\log(S_n^{(\al)}/n)}{n(n+1)(n+2)}-\f2{(n+1)^2}+\f2{(n+2)^2}+\f4{(n+2)^3}+\f{2q_n^{(\al)}}{(n+1)(n+2)}+\f{2(q_n^{(\al)})^2}{n+2}.
\endalign$$
Combining this with (3.2) and (3.3) and noting that $(350001/350000)n^2\gs n(n+1)$, we obtain
$$\align D_n^{(\al)}<&\f{-2\al\log n}{n(n+1)(n+2)}-\f{2(2n+3)}{(n+1)^2(n+2)^2}+\f4{(n+2)^3}
\\&+\f{2c_{\al}}{n(n+1)(n+2)}+\f{2c_{\al}^2}{n^2(n+2)}
\\<&\f{-2\al\log n}{n(n+1)(n+2)}-\f4{(n+1)(n+2)^2}+\f4{(n+1)(n+2)^2}
\\&+\f{2c_{\al}+2(350001/350000)c_{\al}^2}{n(n+1)(n+2)}
\\=&\f{2((350001/350000)c_{\al}^2+c_{\al}-\al\log n)}{n(n+1)(n+2)}.
\endalign$$
Note that
$$\align&\f{350001}{350000}c_{\al}^2+c_{\al}
\\=&\f{350001}{350000}\times1.085^2(\al+1)^2 1.2^{2\al}+1.085(\al+1)1.2^{\al}
\\<&1.2(\al+1)^2 1.2^{2\al}+1.2(\al+1)1.2^{\al}\ls\al\log N(\al)\ls \al\log n.
\endalign$$
So we have $D_n^{(\al)}<0$ and hence
$$\f{\root{n+1}\of{S_{n+1}^{(\al)}/(n+1)}}{\root n\of{S_n^{(\al)}/n}}<\f{\root{n+2}\of{S_{n+2}^{(\al)}/(n+2)}}{\root {n+1}\of{S_{n+1}^{(\al)}/(n+1)}}$$
as desired. \qed

\Ack. The work was done during the author's visit to the University of Illinois at Urbana-Champaign, so the author wishes to thank Prof. Bruce Berndt
for his kind invitation and hospitality.


\widestnumber\key{CP}

 \Refs

\ref\key AS\by M. Abramowitz and I. A. Stegun (eds.)\book Handbook of Mathematical Functions with Formulas, Graphs, and Mathematical Tables
\publ 9th printing, New York, Dover, 1972\endref

\ref\key C\by W. Y. C. Chen\paper Recent developments on log-concavity and $q$-log-concavity of combinatorial polynomials
\jour a talk given at the 22nd Inter. Confer. on Formal Power Series and Algebraic Combin. (San Francisco, 2010)\endref

\ref\key D\by P. Dusart\paper The $k$th prime is greater than $k(\log k+\log\log k-1)$ for $k\gs2$
\jour Math. Comp. \vol 68\yr 1999\pages 411--415\endref

\ref\key HR\by G. H. Hardy and S. Ramanujan\paper Asymptotic formulae in combinatorial analysis
\jour Proc. London Math. Soc.\vol 17\yr 1918\pages 75--115\endref

\ref\key J\by J. E. Janoski\paper A collection of problems in combinatorics\jour PhD Thesis, Clemson Univ., May 2012\endref

\ref\key LW\by L. L. Liu and Y. Wang\paper On the log-convexity of combinatorial sequences\jour Adv. in Appl. Math.\vol 39\yr 2007\pages 453--476
\endref

\ref\key R\by P. Ribenboim\book The Little Book of Bigger Primes\publ 2nd Edition, Springer, New York, 2004\endref

\ref\key Ro\by J. B. Rosser\paper The $n$-th prime is greater than $n\log n$\jour Proc. London Math. Soc. \vol 45\yr 1939\pages 21--44\endref

\ref\key RS\by J. B. Rosser and L. Schoenfeld\paper Approximate formulas for some functions of prime numbers
\jour Illinois J. Math. \vol 6\yr 1962\pages 64--94\endref

\ref\key S\by Z.-W. Sun\paper On functions taking only prime values\jour preprint (2012), {\tt arXiv:1202.6589}\endref


\endRefs
\enddocument